\documentclass[reqno]{amsart}

\RequirePackage[left,mathlines]{lineno}%
\RequirePackage{amsthm}%
\RequirePackage{amsmath}%
\RequirePackage{amsxtra}%
\RequirePackage{amstext}%
\RequirePackage{amssymb}%
\RequirePackage{latexsym}%
\RequirePackage{dsfont}%

\newtheorem{theorem}{Theorem}[section]

\theoremstyle{definition}
\newtheorem{definition}[theorem]{Definition}

\newtheorem{example}[theorem]{Example}

\theoremstyle{remark}

\theoremstyle{corollary}

\numberwithin{equation}{section}

\ifx\pdfoutput\undefined%
   \IfFileExists{graphicx.sty}{\RequirePackage[dvips]{graphicx}
   \DeclareGraphicsExtensions{.eps,.ps}}{}
\else%
   \ifnum\pdfoutput=0
      \IfFileExists{graphicx.sty}{\RequirePackage[dvips]{graphicx}
      \DeclareGraphicsExtensions{.eps,.ps}}{}
   \else
      \IfFileExists{graphicx.sty}{\RequirePackage[pdftex]{graphicx}
      \DeclareGraphicsExtensions{.pdf,.png,.jpg}}{}
   \fi
\fi 

\graphicspath{{./figs/}}%

\RequirePackage[numbers,sort]{natbib}

\RequirePackage[colorlinks,citecolor=blue,urlcolor=blue]{hyperref}
\RequirePackage{hypernat}%
\RequirePackage{color}%

\definecolor{red}       {rgb}{0.0,0.0,1.0}        
\definecolor{magenta}   {rgb}{0.0,0.0,1.0}        
\definecolor{cyan}      {rgb}{0.0,0.0,1.0}        
\definecolor{green}     {rgb}{0.0,0.4,0.3}        

\begin{document}



\title[Control on the second order derivatives]%
{Quest for the control on the second order derivatives: topology
optimization with functional includes the state's curvature}


\author{R. Tavakoli}

\address{Rouhollah Tavakoli, Department of Material Science and Engineering, Sharif University of
Technology, Tehran, Iran, P.O. Box 11365-9466}

\thanks{Rouhollah Tavakoli,
Department of Material Science and Engineering, Sharif University of
Technology, Tehran, Iran, P.O. Box 11365-9466,
\href{mailto:rtavakoli@alum.sharif.edu}{rtavakoli@alum.sharif.edu},
\href{mailto:rohtav@gmail.com}{rohtav@gmail.com},
\href{mailto:tav@mehr.sharif.edu}{tav@mehr.sharif.edu}.}

\email{\href{mailto:rtavakoli@alum.sharif.edu}{rtavakoli@alum.sharif.edu},
\href{mailto:rohtav@gmail.com}{rohtav@gmail.com},
\href{mailto:tav@mehr.sharif.edu}{tav@mehr.sharif.edu}}

\date{\today}%


\maketitle%


\begin{abstract}%

Many physical phenomena, governed by partial differential equations
(PDEs), are second order in nature. This makes sense to pose the
control on the second order derivatives of the field solution, in
addition to zero and first order ones, to consistently control the
underlaying process. However, this type of control is nontrivial and
to the best of our knowledge there is nigher a theoretic nor a
numeric work in this regard. The present work goals to do the first
quest in this regard, examining a problem of this type using a
numerical simulation.

A distributed parameter identification problem includes the control
 on the diffusion coefficient of the Poisson equation and a functional includes
the state's curvature is considered. A heuristic regularization tool
is introduced to manage codimension-one singularities during the
functional analysis. Based on the duality principles, the
approximate necessary optimality conditions is found. The system of
optimality conditions is solved using a globalized projected
gradient method. Numerical results, in two- and three-dimensions,
implied the possibility of posing control on the second order
derivatives and success of the presented numerical method.
\vspace*{2mm}


\end{abstract}%

\tableofcontents%


\section{Introduction}

Many problems in engineering sciences and physics are modeled by
partial differential equations (PDEs). The numerical solution of
PDEs provides useful information for design engineers to predict the
behaviors of corresponding systems. In the highly competitive world
of today, it is no longer sufficient to design a system that
performs the required task satisfactorily, but it is desired to find
 the optimal conditions. The PDE-constrained optimization
 (cf. \cite{lions1971ocs,bendsoee2004tot}) is an
 effective way to approach this goal.

In the context of PDE-constrained optimal control it is common to
pose either pointwise or global control on the state solution and/or
state's derivatives. Most of works in this regard are focused to
apply control on the the zero and/or first order derivatives of
filed solution.

However, most of the physical phenomena are spatial second order in
nature. This means that the behavior of a point is a function of its
neighborhood. This make sense to pose control on the second order
derivatives of the state solution, in addition to the zero and/or
first order derivatives. Moreover, in some cases we may not have
sufficient information about the desired conditions in terms of the
zero and/or first order information. Lets to give an instance, to
make the mentioned problem more clear. In heat transfer problems,
hot-centers (the locus of heat flux concentration) are a subset of
temperature field critical points, i.e., where the temperature
gradient approaches to zero. Now assume we want to control the
position of hot centers. It is clear that to pose the control on the
state's gradient does not make a physically consistent design
problem, because, the cold-centers and saddle points are also
included within the set of mentioned critical points. An effective
way to consistently contrast between these three types of points, is
to consider the second order information too, i.e., to include the
Hessian of the state solution. It is worth mentioning that, some of
the currently used design problems could be enriched/regularized by
including the second order information in their formulation.

However, to pose the control on the second order information is not
an easy task. To the best of our knowledge there is neither a
theoretic nor a numeric work in this regard. Since theoretical study
on this topic is not easy using available functional analysis tools,
we shall perform a numerical quest on this topic here. To be more
specific, we consider a topology optimization problem (cf.
\cite{bendsoee2004tot}) in this study. Our problem includes the
minimization of an integral functional includes the state's first
and second order derivatives where the sate is governed by a Poisson
equation. Moreover the diffusion coefficient of the Poisson equation
is the control parameter in this study. In fact we deal with a
PDE-constrained distributed parameter identification problem in this
study.


\section{Problem statement}%
\label{sec:primal}%

Consider a simply connected domain $Q \subset \mathbb{R}^d$ ($d=$ 2
or 3), with sufficiently regular boundary $\Sigma = \partial Q$,
where the optimization problem will take place. Consider the
following design problem in steady-state heat (or electrical)
conduction: given two isotropic conducting materials, with thermal
conductivities $\alpha$ and $\beta$, $0<k_\beta<k_\alpha$. Suppose
that at each spatial point $\textbf{x} \in Q$ the conductivity is
given in terms of the $\alpha-$phase characteristic function $\chi =
\chi(\textbf{x})$, where $\chi \in \{0, 1\}$, as follows:%
\begin{equation}%
\label{eq:conductivity}%
    k (\chi) = \chi k_\alpha + (1-\chi) k_\beta
\end{equation}%
where $k_\alpha$ and $k_\beta$ denote the thermal (or electrical)
conductivity of phases $\alpha$ and $\beta$ respectively. In fact,
$\chi(\textbf{x})$ plays the role control or design parameter here.
The goal is to find $\chi(\textbf{x})$ configuration which solve the
following design problem: %
\begin{equation}%
\label{eq:P}%
\nonumber%
\mathtt{(P)}:= \ \ \arg \min_{\chi\in \Upsilon}  F(\chi) =
\int_{\mathcal{D}}
 f (\chi, \nabla u, \nabla \nabla u) \ d\Omega
\end{equation}%
subject to:
\begin{equation}%
\label{eq:DP}%
\nabla \cdot \big(k(\chi) \nabla u \big) = g(\mathbf{x}) \ \
\mathtt{in}\ Q, \ \ u = u_0 \ \  \textrm{on} \ \Sigma   %
\end{equation}%

where $\mathcal{D} \subseteq \Omega$ denotes the support domain of
the integral functional, $d \Omega$ denotes the volume measure
induced in $Q$, $u \in X_u(Q)$ is the state function, $\nabla u$ is
the state gradient, $\nabla \nabla u$ denotes the local Hessian
matrix of the state function, i.e., $\nabla \nabla u =
[u_{ij}]_{d\times d}$, $i,j \in\{x,y\}$ ($u_{ab}= \frac{\partial^2
u}{\partial a \partial b}$), $g \in X_g(Q)$ is a given source
function, $u_0(\textbf{x}) \in X_u(Q)$ denotes the Dirichlet
boundary condition and $\Upsilon$ denotes the admissible design
 space which is defined as follows%
\begin{equation}%
\label{eq:adspace}%
    \Upsilon := \big\{  \chi \in X_\chi(Q) \
    \big| \quad
    \ R_L |Q| \leq \int_Q \chi d\Omega \leq R_U |Q|,
    \quad \chi \in \{0. 1\}\  \big\}
\end{equation}%
where $0<R_L\le R_U<1$ are given lower and upper bounds on
$\alpha$-phase material resource respectively, and $|\Omega|$
denotes the measure of $\Omega$ (volume for $d=3$ and area for
$d=2$). In above formulae $X_u$, $X_g$ and $X_\chi$ denotes some
appropriate Banach spaces. Since our goal in this study is not to do
a rigorous mathematical analysis, we simply assume that these
function spaces possess sufficiently regularity required by our
analysis. It is well known that optimal design problems like
$\mathtt{(P)}$, are ill-posed and do not admit an optimal solution
in class of characteristic functions (cf. \cite{ekeland1976convex}).
The relaxation form of these problem is usually achieved extending
the design domain $\Xi$ to the convex and continuous space $\Xi$
which is defined as follows,
\begin{equation}%
\label{eq:radspace}%
    \Xi =
    \big\{  w \in X_w(Q) \
    \big| \quad
    \ R_L |Q| \leq \int_\Omega w \ d \textbf{x} \leq R_U |Q|,
    \quad 0 \le w \le 1  \big\}
\end{equation}%
where $X_w$ is a sufficiently regular Banach space. In fact in the
relaxed design problem the material conductivity could varies
gradually between phases $\alpha$ and $\beta$ and we have somehow a
functionally graded material. The SIMP approach
\cite{bendsoee2004tot} is employed in this study to achieve a near
0-1 topology. In this method the conductivity function
\eqref{eq:conductivity} is replaced by the following function:%
\begin{equation}%
\label{eq:ksimp}%
    k(w) = w^q (k_\alpha - k_\beta) + k_\beta
\end{equation}%
where $q \geq 1$ is a penalization factor which is commonly called
as the SIMP power.

In the present study, the function $f$ in $\mathtt{(P)}$ is assumed
to have the following form: %
\[f := a w^b (\kappa-\kappa_0)^c,\]%
where $a, b, c \in \mathbb{R}$, $b, c \ge 1$ are some given
constants, $\kappa$ is the local mean curvature of state
iso-contours, i.e., $\kappa = \nabla \cdot (\nabla u / |\nabla u|)$.
It is clear that $\kappa$ is a nonlinear function of both the first
and second order derivatives of $u$ (cf. \cite{osher2003lsm}).
Therefore, this type of functional it is a fairly good choice for
the purpose of this study. Besides it has a clear physical
interpretation and has certain applications in practical design
problems. For example, a negative curvature with large absolute
value is signature of thermal centers in the heat transfer problems
(cf. \cite{tavakoli2010prediction}).


\section{Regularization of the curvature singularities}%
\label{sec:kappareg}%

Assume $u_0$ is a constant function on $\Gamma$, then according to
the maximum principles, there is at least a point $\mathbf{x}_0 \in
\Omega$ such that $|\nabla u| = 0$). To avoid the curvature
singularity in this regard, we simply replace $|\nabla u|$ by
$|\nabla u|_\varepsilon := \sqrt{\varepsilon^2 + |\nabla u|^2}$ in
the definition of mean curvature. Therefore the regularized mean
curvature is defined as follows: $\kappa_\varepsilon = \nabla \cdot
(\nabla u / |\nabla u|_\varepsilon)$, where the small number
$\varepsilon \in \mathbb{R}^+$ is the curvature regularization
parameter. For the purpose of convenience we simply use $\kappa$
instead of $\kappa_\varepsilon$ henceforth. In practice we use a
sufficiently small value ($10^{-20}$ in this study) instead of
$\varepsilon$ without any  outer iteration on $\varepsilon$ to
recover the original problem.

Consider the discretized version of the optimal design problem with
the minimum grid size $h$. According to \cite{osher2003lsm}, the
maximum absolute value of the curvature which can be resolved by
such a grid resolution is equal to $\frac 1 h$. Therefore, if the
computed value of $\kappa$ locates outside of $-\frac 1 h \le \kappa
\le \frac 1 h$, we merely replace that value with either $-\frac 1
h$ or $\frac 1 h$ depending on the sign of curvature.


\section{A Local analysis near a shallow gradient regions}%
\label{sec:la}%

Suppose that $u$ is $C^2$ near a generic spatial stationary point
$\textbf{x}_0$, i.e., $|\nabla u| = 0$ and $\nabla \nabla u$ is
nonsingular there. Consider d-dimensional closed ball $\mathcal{B}_r
\subset Q$ with radius $r \in \mathbb{R}^+$ centered at
$\mathbf{x}_0$. We would like to study the behavior of the objective
functional in $(\mathtt{P})$ restricted to ball $\mathcal{B}_r$. In
this regard we can restrict ours to case $\kappa_0=0$ i.e., to study
the behavior of the following functional:
\begin{equation}%
\label{eq:RVPregB}%
    F |_{\mathcal{B}_r} :=
    \int_{\mathcal{B}_r}
    a w^b \kappa^c \ d\Omega %
\end{equation}
Since $\kappa$ is invariant under an Euclidean transformation, we
assume $\mathbf{x}_0 = \mathbf{0}$. Consider an approximation of $u$
in the neighborhood of $\mathbf{x}_0$ by the Taylor
expansion up to the second order derivatives (recall that $|\nabla u(\mathbf{x}_0)| = 0$):%
\begin{equation}%
\label{eq:taylorseris}%
    u(\textbf{x}) \approx
    u(\mathbf{x}_0) +
    \frac{1}{2}\ \textbf{x}^T \cdot
    \nabla \nabla u (\mathbf{x}_0) \cdot \textbf{x}%
\end{equation}%
Since $\kappa$ is rotationally invariant, we can assume that
\begin{equation}%
\label{eq:diaglambda}%
   \nabla \nabla u(\mathbf{x}_0) =
   \texttt{diag}(\lambda_0(\mathbf{x}_0), \cdots, \lambda_d(\mathbf{x}_0))
\end{equation}%
where $\lambda_i(\mathbf{x}_0)$ denotes the $i$-th eigenvalue of
$\nabla \nabla u(\mathbf{x})$ at $\mathbf{x} = \mathbf{x}_0$.
Therefore, by \ref{eq:taylorseris}
we have%
\begin{equation}%
\label{eq:taylorseris2}%
    u(\textbf{x}) \approx
    u(\mathbf{x}_0) +
    \frac{1}{2}\ \sum_{i=1}^d \lambda_i(\mathbf{x}_0) x_i^2 %
\end{equation}%
The local structure of the solution in the vicinity of
$\textbf{x}_0$ is a function of relative values of
$\lambda_i(\mathbf{x}_0)$. There are four generic types of local
structures in this regard: sheet-like, tube-like, blob-like,
double-cone-like local structures. Among these local features, the
blob (ellipsoidal) and double-cone (hyperbolic) types have the
maximum mean curvature concentration. For the ease of analysis, we
restrict ours to the blob-like local structures (when $\lambda_1
\approx \cdots \approx \lambda_d$). Assume $\lambda_1 \approx \cdots
\approx \lambda_d = \lambda$, i.e., a spherical local structure,
using \eqref{eq:taylorseris2} and straightforward computations, we
have:
\begin{equation}%
\label{eq:kappaball}%
   \kappa = \frac{d-1}{\sqrt{\sum_{i=0}^d x_i^2}}
\end{equation}
substitution of \ref{eq:kappaball} in \ref{eq:RVPregB} results:%
\begin{equation}%
\label{eq:RVPregB2}%
    F |_{\mathcal{B}_r} =
    \int_{\mathcal{B}_r}
    a w^b \bigg(\frac{d-1}{\sqrt{\sum_{i=0}^d x_i^2}}\bigg)^c \ d\Omega %
\end{equation}
Considering the symmetric structure of the blob-like structure and
using the spherical coordinate results:
\begin{equation}%
\label{eq:RVPregB3}%
    F |_{\mathcal{B}_r} =
    \int_{0}^r
    a w^b \bigg(\frac{d-1}{r}\bigg)^c \ \pi (2 r)^{d-1} dr %
    = \zeta \int_{0}^r r^{d-c-1} \ dr
\end{equation}
where constant $\zeta$ is equal to $2^{(d-1)} (d-1)^c \pi a w^b$.
Considering \ref{eq:RVPregB3}, integral functional \ref{eq:RVPregB}
is well-defined in the vicinity of a regular shallow gradient point
if $c \le d-1$. It turns out that we are not free to select the
power of the curvature term in our integral functional. For
instance, assume we are interested to solve a topology optimization
problems subject to pointwise constraints on mean curvature.
According to above analysis, in the case of two dimensions, we can
not use a quadratic (or higher order) penalty function (cf.
\cite[][ch. 3]{bertsekas1982coa}) to manage local constraints (it is
worth mentioning that using quadratic penalty functions is very
common to manage such problems. In the other word we have to use
exact penalty functions (cf. \cite[][ch. 4]{bertsekas1982coa}) in
these cases.


\section{Regularization of singular integrals}%
\label{sec:regsi}%

To derive the first order necessary optimality conditions we have to
deal with some singular integrals to continue the analysis. The good
of this section is to introduce some (heuristic) regularization
tools to manage such cases. Let $f(\textbf{x}): \mathbb{R}^d
\rightarrow \mathbb{R}$ be a real-valued function defined on the
codimension-one $C^2$ manifold $\mathcal{M} \subset \mathbb{R}^d$
which is embedded within domain $Q$. We are interested to compute
the following singular
integral:\\%
\begin{equation}%
\label{eq:si}%
   I := \int_{\mathcal{M}} f(\textbf{x}) \ d \Gamma = %
       \int_{Q} f(\textbf{x}) \ \delta (\mathcal{M}) \ d \Omega  %
\end{equation}
where $d\Gamma$ denotes the surface measure induced on
$\mathcal{M}$, $\delta (\mathcal{M})$ denotes the delta-distribution
concentrated on $\mathcal{M}$. Consider $b_{\mathcal{M}} \in
C^1(\Omega)$ as the Euclidean signed distance function with respect
to manifold $\mathcal{M}$. Assume $(\cdot)^{ext}_\sigma$ denotes the
smooth extension (along the normal direction) of scaler field
$(\cdot)$ defined on $\mathcal{M}$ into
$\mathcal{N}(\mathcal{M},\sigma)$, where
$\mathcal{N}(\mathcal{M},\sigma)$ denotes $\pm \sigma$-distance
neighborhood of $\mathcal{M}$. Moreover assume that the manifold
forms the boundaries of $\mathcal{N}(\mathcal{M},\sigma)$ is
non-degenerate (there is no self-crossing). Assume that
$\mathcal{N}(\mathcal{M},\sigma)$ is spanned by $s\lambda$-system
such that $s$ coincides to the iso-contours of the distance function
$b_{\mathcal{M}}$ and $\lambda$ coincides to normals on the
$b_{\mathcal{M}}$-constant surfaces. Therefore, $ds$ is a
$(d-1)$-dimensional (area) metric and $d\lambda$ is a
one-dimensional metric. It is obvious that $s\lambda$-system
consistently spans $\mathcal{N}(\mathcal{M},\sigma)$, if
 iso-contours of $b_{\mathcal{M}}$ do not cross each other (which is
followed by our assumption). Therefore we can use $s\lambda$-system
together with its corresponding d-dimensional volume element,
$dsd\lambda$, to compute $d$-dimensional volume integrals. Applying
 $s\lambda$-system to integral \eqref{eq:si} results:%
\begin{equation}%
\label{eq:casi}%
   I_\sigma = %
       \int_{\lambda=-\sigma}^{\lambda=+\sigma}%
       \rho_\sigma(\lambda)%
       \int_{\mathcal{M}_{b} : \ b_{\mathcal{M}} = \lambda}%
       f^{ext}_\sigma(\textbf{x})%
       \ ds\ d\lambda%
\end{equation}
\\ where $\lim_{\sigma \rightarrow 0}
I_\sigma = I$, the manifold $\mathcal{M}_{b} \subset
\mathcal{N}(\mathcal{M},\sigma)$ coincides to an $s$-surface which
is located at distance $b$ from $\mathcal{M}$, and the weighting
function $\rho_\sigma \in C^\infty(\mathbb{R})$ is chosen such that
the intensity (measure) of $f(\textbf{x})$ be preserved, i.e.,
\begin{equation}%
\label{eq:kernelprop1}%
       \int_{\lambda=-\sigma}^{\lambda=+\sigma}%
       \rho_\sigma(\lambda) \ d\lambda = 1%
\end{equation}%
the other requirements of function $\rho_\sigma$ are:
$\rho_\sigma(x) \neq 0$ for $x \in (-\sigma, +\sigma)$,
$\rho_\sigma(x) = 0$ for $x \notin (-\sigma, +\sigma)$. In fact,
$\rho_\sigma$ makes a converging sequence to the classic delta
function, $\delta$, in the sense that: $\lim_{\sigma \rightarrow 0}
\rho_\sigma(x) = \delta(x)$ (cf. \cite{boykin2003ddd}). in the
present study, the following normalized Gaussian distribution is
used as $\rho_\sigma$ function:%
\begin{equation}%
\label{eq:mollifier}%
    \rho_\sigma(\lambda)= \newline%
    \left\{%
\begin{array}{lll}
      \rho_\sigma^0\ \sigma^{-d}\
      \exp \big(  \frac{\sigma^2}{|\lambda|^2 - \sigma^2} \big),
      & \texttt{if} & |\lambda| < \sigma \\
      0,& \texttt{if} & |\lambda| \ge \sigma
\end{array}%
\right.%
\end{equation}%
where $\rho_\sigma^0$ is determined bases on constraint
\eqref{eq:kernelprop1}. Note that $\rho_\sigma$ plays also the role
of a regularization kernel in our formulation (cf.
\cite{temple1955tgf,colombeau1990md}). More precisely, for every
point $\textbf{x}_0 \in \mathcal{M}$, $\rho_\sigma$ does convolution
of $f(\textbf{x}_0)$ along the line segment $(-\sigma,+\sigma)$
which is normal to $\mathcal{M}$ and passes from $\textbf{x}_0$.
Therefore, if the extension of $f$ to neighborhood of the desired
manifold be not smooth, e.g., there are some equi-distanced points
from the manifold surfaces, our method performs well.

Now, lets to back from $s\lambda$-coordinate to our original
Cartesian coordinate. The coarea formula \cite{fleming1960ift} is
used for this purpose. Assume that the change of variable is global
in $\Omega$, the coarea formula can be
expressed as:%
\begin{equation}%
\label{eq:coarea}%
   \int_\Omega |\nabla b_{\mathcal{M}}| \ d\Omega = \int \ ds\
   d\lambda%
\end{equation}%
i.e., $|\nabla b_{\mathcal{M}}| \ d\Omega = ds\ d\lambda$ (cf.:
\cite[][ch. 3]{evans1992mta}). Applying \eqref{eq:coarea} to
\eqref{eq:casi} results:%
\begin{equation}%
\label{eq:casi2}%
   I_\sigma = %
       \int_{\mathcal{N}(\mathcal{M},\sigma)}
       \rho_\sigma(b_{\mathcal{M}}(\textbf{x}))\ %
       f^{ext}_\sigma(\textbf{x})\ %
       |\nabla b_{\mathcal{M}}|%
       \ d \Omega %
\end{equation}
since $\rho_\sigma =0$ for $ \mathbf{x} \in \Omega \setminus
\mathcal{N}(\mathcal{M},\sigma)$,%
\begin{equation}%
\label{eq:casi3}%
   I_\sigma = %
       \int_\Omega
       \rho_\sigma(b_{\mathcal{M}}(\textbf{x}))\ %
       f^{ext}_\sigma(\textbf{x})\ %
       |\nabla b_{\mathcal{M}}|%
       \ d \Omega %
\end{equation}
it is worth noticing that a formula similar to \eqref{eq:casi3} is
derived heuristically (not rigorously) in \cite[][ch.
1]{osher2003lsm}.

The Euclidean signed distance function $b_{\mathcal{M}}$ can be
efficiently computed solving the following eikonal equation
\cite{sethian1999fmm,osher2003lsm}:%
\begin{equation}%
\label{eq:eikonal}%
 |\nabla b_{\mathcal{M}} | =  1 \ \  \mathtt{in}\ \  Q, %
 \quad b_{\mathcal{M}} = 0 \ \  \mathtt{on} \ \   \mathcal{M}
\end{equation}%
this equation can be solved efficiently using either the fast
marching method \cite{sethian1999fmm}. Notice that the solution of
\eqref{eq:eikonal} is a member of $BV(Q)$, i.e., it is not
essentially sufficiently smooth. However, due to the application of
 the regularization kernel $\rho_\sigma$ in our formulation, we are not
concerned about the smoothness of $b_{\mathcal{M}}$. We also need to
extend function $f$ defined on $\mathcal{M}$ along the normal
direction into $\mathcal{N}(\mathcal{M},\sigma)$. The following
PDE is suggested in \cite{adalsteinsson1999fce} to do this job:%
\begin{equation}%
\label{eq:fextension}%
f^{ext} \cdot \nabla b_{\mathcal{M}}  =  0  \ \  \mathtt{in} \ \ Q,
\ \ f^{ext} = f \ \ \mathtt{on} \ \  \mathcal{M}
\end{equation}%
where $f^{ext}$ denotes the extension of $f$ inside
$\mathcal{N}(\mathcal{M},\sigma)$. In \cite{adalsteinsson1999fce},
an efficient fast marching method is suggested to compute
$b_{\mathcal{M}}$ and $f^{ext}$ simultaneously. The implementation
of this method will be used in this study.

For the convenience in notation, the following notation is used in
this study, everywhere required:%
\begin{equation}%
\label{eq:regdef}%
   \int_{\mathcal{M}} (\cdot) \ d \Gamma%
    \equiv
   \int_{\mathcal{N}(\mathcal{M},\sigma)} \langle\langle\ (\cdot) \ \rangle\rangle   \ d \Omega
\end{equation}%
Let $a(\textbf{x})$ and $b(\textbf{x}):\mathbb{R}^d \rightarrow
\mathbb{R}$ and $\textbf{c}(\textbf{x})$ and
$\textbf{d}(\textbf{x}): \mathbb{R}^d \rightarrow \mathbb{R}^d$ be
functions defined on co-dimension one manifold $\mathcal{M}$. To do
our sensitivity analysis, we need the following properties for
operator $\big<\big<(\cdot)\big>\big>$ (note that these relations
are not exactly hold however we are hopeful, without proof, to
validity of them in the sense of  approximation):%
\begin{eqnarray}%
\label{eq:regprop7}%
\label{iden:1}
   \int_{(\cdot)} \langle\langle\ a \ b \ \rangle\rangle\ d \Omega &\approx&%
   \int_{(\cdot)} \langle\langle\ a \ \rangle\rangle \ b \ d \Omega
    \\ \label{iden:2}%
   \int_{(\cdot)} \langle\langle \ \textbf{c} \cdot \textbf{d}\ \rangle\rangle \ d \Omega &\approx&%
   \int_{(\cdot)} \langle\langle \ \textbf{c}\ \rangle\rangle \cdot \textbf{d} \ d \Omega
\end{eqnarray}%
%
%


\section{The first order necessary optimality conditions}%
\label{sec:sensitivity}%

There are several methods to solve a distributed parameter
identification problem. Due to large number of control parameters, a
gradient based method is adapted in this study (cf.
\cite{allaire2007numerical}). Lets to define the meaning of
differentiability on function spaces. There exist several
differentiability notions in mathematical programming literature.
The notion of G{\^a}teaux derivative is applied here.

\begin{definition}%
\label{def:gato}%
(G{\^a}teaux derivative, cf. \cite{allaire2007numerical}) Consider
Banach spaces $Y$ and $W$ and $U$ as open subset of $Y$. A function
$f: U \rightarrow W$ is called to be G{\^a}teaux differentiable at
$u \in U$ if for every test function $v \in Y$ the following limit exist: %
\[
f^\prime(u) :=  \lim_{\zeta \rightarrow 0} \frac{f(u+ \zeta v) -
f(u)}{\zeta}
\]
In this case we show the G{\^a}teaux  derivative symbolically by
$f^\prime(u)$. If $Y$ is a Hilbert space, which is assumed in this
study, then $f^\prime(u)$ lives on the dual space of $Y$. Therefore
using the Riesz representation theorem, there is a unique $e \in Y$
such that $\langle e, v \rangle = f^\prime(u)$, where $\langle
\cdot, \cdot \rangle$ denotes the inner product on $Y$. In this
case, without confusion, it is common to call $e$ as the G{\^a}teaux
derivative. We use notation $d f(u)$ in this case, i.e., $d f(u) =
e$. It is easy to verify that under some mild conditions (which
usually hold in practice) most of properties for classical
derivatives have equivalent extension to G{\^a}teaux derivative. In
the case of multi-variable functional, the partial G{\^a}teaux
derivative are denoted by $f^\prime_{(\cdot)}$, $\partial_{(\cdot)}
f(\cdots)$ symbols in this study. For the purpose of convenience the
G{\^a}teaux derivative is called as the directional derivative in
this study, henceforth. Without confusion, the notation $\langle
\cdot, \cdot \rangle_Q := \int_Q (\cdot) (\cdot) \ d \Omega$ also
used to denote the duality pairing on function spaces in this study. %
\end{definition}

Consider an arbitrary function $p \in X_p(Q)$, where $X_p$ is a
sufficiently regular Banach space. Lets to introduce the following
lagrangian augmenting the inner product of $p$ and the Poisson
equation to the objective functional in problem $\mathtt{(P)}$:
\[
\mathcal{L}(w, u, p) := F(w, u) + %
%
\big\langle \nabla \cdot (k \nabla u)- g, \ p  \big\rangle_Q +
\big\langle u-u_0, \ p  \big\rangle_\Sigma
\]
The set of points satisfy the first order necessary optimality
conditions of problem $\mathtt{(P)}$, denoted by $\mathcal{O}$ can
be expressed as follows (cf. \cite{allaire2007numerical}): %
\[
    \mathcal{O} := %
    \left\{%
    (w,u,p) \in \big(\Xi \times X_u(Q) \times X_p(Q) \big) \
    \Bigg|
\begin{array}{r}
\partial_w \mathcal{L}(w,u,p)=0 \ \ \texttt{in} \  Q  \ \ (\mathtt{C.1})\\%
\partial_u \mathcal{L}(w,u,p)= 0 \ \ \texttt{in} \ Q  \ \ (\mathtt{C.2})\\%
\partial_p \mathcal{L}(w,u,p)=0 \ \ \texttt{in} \  Q  \ \  (\mathtt{C.3})%
\end{array}%
\right\}.%
\]
where $\partial_w \mathcal{L}$, $\partial_u \mathcal{L}$ and
$\partial_p \mathcal{L}$ denote respectively the partial directional
derivatives of $\mathcal{L}$ with respect to $w$, $u$ and $p$ along
arbitrary test directions $\delta w \in X_w(Q)$, $\delta u \in
X_u(Q)$ and $\delta p \in X_p(Q)$ respectively. In fact,
$\mathcal{O}$ includes constrained stationary points of augmented
lagrangian $\mathcal{L}$. Regarding to $\partial_w \mathcal{L}$
in condition $(\mathtt{C.1})$ of set $\mathcal{O}$ we have, %
\begin{align}%
\label{eq:dw_lag}%
\big\langle \partial_w \mathcal{L}, \delta w \big\rangle_Q & =
\big\langle  \partial_w f, \ \delta w \big\rangle_{\mathcal{D}} +  %
\big\langle p\nabla\cdot(\partial_w k\nabla u),\ \delta w
\big\rangle_Q := I_1 + I_2
\end{align}%
where $\partial_w f = a b w^{(b-1)} (\kappa-\kappa_0)^c$ and
$\partial_w k = q w^{(q-1)} (k_\alpha - k_\beta)$. For $I_2$ we
have,
\begin{align}%
\label{eq:dw_lag_I2}%
I_2 =
\big\langle p \partial_w k \nabla u \cdot \mathbf{n},\ \delta w
\big\rangle_\Sigma
-
\big\langle \partial_w k \nabla u \cdot \nabla p,\ \delta w
\big\rangle_Q
:= I_3 + I_4
\end{align}%
where $\mathbf{n}$ denotes the outer unit normal on $\Sigma$.

Regarding to $\partial_u \mathcal{L}$ in condition $(\mathtt{C.2})$
of set $\mathcal{O}$ we have, %
\begin{align}%
\label{eq:du_lag}%
\big\langle \partial_u \mathcal{L}, \delta u \big\rangle_Q & =
\big\langle  \partial_u f, \ \delta u \big\rangle_{\mathcal{D}} +  %
\big\langle \nabla\cdot(k\nabla \delta u),\ p \big\rangle_Q +
\big\langle \delta u,\ p \big\rangle_\Sigma
:= I_5 + I_6 + I_7
\end{align}%
Applying two consecutive times  the integration by part followed by
the divergence theorem to term $I_6$ in \eqref{eq:du_lag} results:%
\begin{align}%
\label{eq:du_lag_I6}%
I_6 =  %
\big\langle k \nabla\delta u \cdot \mathbf{n}, \  p
\big\rangle_\Sigma -
\big\langle k \nabla p \cdot \mathbf{n}, \  \delta u
\big\rangle_\Sigma + %
\big\langle \nabla\cdot(k\nabla p),\ \delta u \big\rangle_Q
:= I_8 + I_9 + I_{10}
\end{align}%
Since $u$ is fixed on $\Sigma$, we have $\delta u = 0$ on $\Sigma$
therefore $I_8 = I_9=0$ in \eqref{eq:du_lag_I6}.
For term $I_5$ in \eqref{eq:du_lag} we have: %
\begin{align}%
\label{eq:du_lag_I5}%
I_5  =  %
\big\langle \partial_\kappa f \partial_u \kappa, \  \delta u
\big\rangle_{\mathcal{D}} =
\int_{\mathcal{D}} \partial_\kappa f \
\nabla \cdot \bigg( \frac{\nabla\delta u}{|\nabla u|} - %
\frac{\nabla u ( \nabla u \cdot \nabla\delta u)} {|\nabla u |^3}
\bigg) \ d\Omega := J_1
\end{align}%
To make expression concise, lets to define the following local
operator $P(u)$, %
\begin{equation}%
\label{eq:projmat}%
P(u) := \mathds{1}_d - \frac{\nabla u}{|\nabla u|}\otimes \frac{\nabla u }{|\nabla u|}%
\end{equation}
where $\mathds{1}_d \in \mathbb{R}^{d\times d}$ denotes the identity
matrix and $\otimes$ denotes the tensor product of $d$-vectors,
i.e., $\textbf{a} \otimes \textbf{b} = [a_i b_j]_{d \times d}$.
Using \eqref{eq:projmat} in \eqref{eq:du_lag_I5}, $J_1$ can be
written as follows:%
\begin{align}%
\label{eq:du_lag_J1}%
J_1 =  %
\int_{\mathcal{D}}%
\partial_\kappa f  \ \nabla \cdot \bigg(\frac{ \nabla\delta u \cdot P(u)}{|\nabla u|}
\bigg) \ d\Omega%
\end{align}%
Applying the integration by part and the divergence theorem to \eqref{eq:du_lag_J1} results: %
\begin{align}%
\label{eq:du_lag_J1_2}%
J_1 =  %
\int_{\partial \mathcal{D}}%
\frac{ \partial_\kappa f \ \nabla\delta u \cdot
P(u)\cdot \mathbf{m}}{|\nabla u|} \ d\Gamma - %
\int_{\mathcal{D}}%
\frac{ \nabla\delta u \cdot P(u) \cdot \nabla (\partial_\kappa
f)}{|\nabla u|}  \ d\Omega  %
:= J_2 - J_3
\end{align}%
where $\mathbf{m}$ denotes the outer unit normal on $\partial
\mathcal{D}$ and $\partial \mathcal{D}$ denotes the sufficiently
regular codimension-one manifold forms the boundaries of
$\mathcal{D}$.
Now lets to proceed the derivation, using regularization tools
introduced in section \ref{sec:regsi}. For $J_2$ in
\eqref{eq:du_lag_J1_2} we have
\begin{align}%
\label{eq:du_lag_J2}%
\nonumber %
J_2 & \stackrel{\eqref{eq:regdef}}{=}  %
\int_{\mathcal{X}}%
\bigg\langle\bigg\langle \frac{ \partial_\kappa f \ \nabla\delta u
\cdot P(u)\cdot \mathbf{m}}{|\nabla u|}
\bigg\rangle\bigg\rangle \ d\Omega %
\stackrel{\eqref{iden:2}}{\approx}
\int_{\mathcal{X}}%
\bigg\langle\bigg\langle\frac{ \partial_\kappa f P(u)\cdot
\mathbf{m}}{|\nabla
u|} \bigg\rangle\bigg\rangle \cdot \nabla\delta u  \ d\Omega \\ \nonumber  %
&=
\int_{\partial \mathcal{X}}%
\delta u\  \bigg\langle\bigg\langle\frac{ \partial_\kappa f
P(u)\cdot \mathbf{m}}{|\nabla
u|} \bigg\rangle\bigg\rangle \cdot \mathbf{k}   \ d\Gamma %
 \ -
\int_{\mathcal{X}}%
  \nabla \cdot \bigg\langle\bigg\langle\frac{ \partial_\kappa f
P(u)\cdot \mathbf{m}}{|\nabla u|} \bigg\rangle\bigg\rangle \delta u   \ d\Omega \\ %
& := J_4 - J_5 %
\end{align}%
where $\mathcal{X} :=\mathcal{N}(\partial \mathcal{D}, \sigma)$ and
$\mathbf{k}$ denotes the outer unit normal on boundaries of
$\mathcal{X}$. For $J_4$ in \eqref{eq:du_lag_J2} we have: %
\begin{align}%
\label{eq:du_lag_J4}%
\nonumber %
J_4 & \stackrel{\eqref{eq:regdef}}{=}  %
\int_{\mathcal{Y}}%
\bigg\langle\bigg\langle\delta u \bigg\langle\bigg\langle\frac{
\partial_\kappa f P(u)\cdot \mathbf{m}}{|\nabla u|}
\bigg\rangle\bigg\rangle \cdot \mathbf{k} \bigg\rangle\bigg\rangle
\ d\Omega \\ %
& \stackrel{\eqref{eq:regdef}}{\approx}  %
\int_{\mathcal{Y}}%
\bigg\langle\bigg\langle  \bigg\langle\bigg\langle\frac{
\partial_\kappa f P(u)\cdot \mathbf{m}}{|\nabla
u|} \bigg\rangle\bigg\rangle \cdot \mathbf{k} \bigg\rangle\bigg\rangle\ \delta u  \ d\Omega %
 := J_6 %
\end{align}%
where $\mathcal{Y} :=\mathcal{N}(\partial \mathcal{X}, \sigma)$.
Applying the integration by part and the divergence theorem to term
$J_3$ in \eqref{eq:du_lag_J1_2} results:%
\begin{align}%
\label{eq:du_lag_J3}%
\nonumber
J_3 & =  %
\int_{\partial\mathcal{D}} \frac{\delta u  \nabla (\partial_\kappa
f) \cdot P(u) \cdot \mathbf{m} }{|\nabla u|}  \ d\Gamma  %
-
\int_{\mathcal{D}}  \nabla \cdot \bigg(\frac{  \nabla
(\partial_\kappa
f) \cdot P(u)}{|\nabla u|}\bigg)\delta u  \ d\Omega := J_7 - J_8 \\  %
& \stackrel{\eqref{eq:regdef}, \eqref{iden:1}}{\approx}  %
\int_{\mathcal{X}} \bigg\langle\bigg\langle\frac{\nabla
(\partial_\kappa
f) \cdot P(u) \cdot \mathbf{m} }{|\nabla u|}\bigg\rangle\bigg\rangle \  \delta u \ d\Omega  %
- J_8 %
\end{align}%
Collecting terms include implicit function $\delta u$, we can solve
the following adjoint Poisson equation to enforce the condition
 $\mathtt{C}.2$ of set $\mathcal{O}$: %
\begin{align}%
\label{eq:ah}%
\nabla \cdot (k(w) \nabla p) = %
h(u) \ \ \mathtt{in} \ \ Q, \quad p = 0 \ \ \mathtt{on} \ \
\Sigma %
\end{align}%
where $u$ solves the direct problem \eqref{eq:DP}, $h(u)
= h_1(u) - h_2(u) - h_3(u) + h_4(u)$ and, %
\begin{align}%
\nonumber%
h_1(u) &:=  \bigg\langle \bigg\langle \ \bigg\langle\bigg\langle
\frac{
\partial_\kappa f P(u)\cdot \mathbf{m}}{|\nabla u|} \bigg\rangle\bigg\rangle
\cdot \mathbf{k} \bigg\rangle\bigg\rangle\ \mathcal{I}_{\mathcal{Y}}(\mathbf{x}) \\ \nonumber%
h_2(u) &:= \nabla \cdot \bigg\langle\bigg\langle\frac{
\partial_\kappa f P(u)\cdot \mathbf{m}}{|\nabla u|} \bigg\rangle\bigg\rangle \
\mathcal{I}_{\mathcal{X}}(\mathbf{x})  \\ \nonumber%
h_3(u) & := \bigg\langle\bigg\langle\frac{\nabla (\partial_\kappa f)
\cdot P(u) \cdot \mathbf{m} }{|\nabla u|}\bigg\rangle\bigg\rangle  \
\mathcal{I}_{\mathcal{X}}(\mathbf{x})  \\ \nonumber%
h_4(u) & := \nabla \cdot \bigg(\frac{  \nabla (\partial_\kappa f)
\cdot P(u)}{|\nabla u|}\bigg)  \
\mathcal{I}_{\mathcal{D}}(\mathbf{x})
\end{align}
where $\mathcal{I}_{\mathcal{D}}(\mathbf{x})$,
$\mathcal{I}_{\mathcal{X}}(\mathbf{x})$ and
$\mathcal{I}_{\mathcal{Y}}(\mathbf{x})$ denote the characteristic
functions of the spatial domains $\mathcal{D}$, $\mathcal{X}$ and
$\mathcal{Y}$ respectively, i.e., for instance
$\mathcal{I}_{\mathcal{D}}(\mathbf{x}) = 1$ for all $\mathbf{x} \in
\mathcal{D}$ and $\mathcal{I}_{\mathcal{D}}(\mathbf{x}) = 0$
elsewhere. According to our numerical experiments, $h_1(\mathbf{x})
\ll h_i(\mathbf{x})$ $(i=2,3,4)$ such that $h_1$ can be ignored
without a sensible loss in the accuracy of numerical solutions.
Since $p=0$ on $\Sigma$, integral $I_3$ in \eqref{eq:dw_lag_I2} is
equal to zero.

It is evident that when $u$ solves the direct Poisson equation
\eqref{eq:DP} the condition $\mathtt{C}.3$ holds in $\mathcal{D}$.
Putting altogether, the (approximate) first order necessary
optimality conditions for optimization problem  $\mathtt{(P)}$ can
be expressed as follows:%
$\mathtt{(OC)}$ as follows:%
\begin{equation}%
\label{eq:opt_cond}%
  \nonumber
    \mathtt{(OC)} := %
    \left\{%
\begin{array}{rll}
    \nabla \cdot (k(w) \nabla u)  = & g(\mathbf{x})  & \texttt{in} \quad Q \\
    u(\mathbf{x})  = & u_0(\mathbf{x}) & \texttt{on} \quad \Sigma \\%
   \nabla \cdot (k(w) \nabla p)  = &  h(\mathbf{x}) & \texttt{in} \quad Q \\
   p(\mathbf{x})  = & 0 & \texttt{on} \quad \Sigma \\ %
     \mathcal{P}_{\Xi}(w - \partial_w \mathcal{L}) - w = & 0 & \texttt{in} \quad
     Q
\end{array}%
\right.%
\end{equation}%
where $\partial_w \mathcal{L}(\mathbf{x}) = j_1(\mathbf{x}) -
j_2(\mathbf{x})$, $j_1(\mathbf{x}) = \partial_w f \
\mathcal{I}_{\mathcal{D}}(\mathbf{x})$, $j_2(\mathbf{x})=
\partial_w k \nabla u \cdot \nabla p$ and operator
$\mathcal{P}_{\Xi}(\cdot)$ denotes the orthogonal projection onto
the admissible control space $\Xi$. Due to the convexity of $\Xi$
the optimality conditions $\mathtt{(OC)}$ is stated in terms of the
projected gradient with respect to the admissible control space. Due
to specific structure of $\Xi$ the projection
$\mathcal{P}_{\Xi}(\cdot)$ can be computed very efficiently in
practice.

It is worth mentioning that in the special case when  $\mathcal{D}:=
Q$, since $\delta u =0$ on $\Sigma$ we have $J_2=J_7=0$. Therefor we
can solve the adjoint Poisson equation exactly without any
regularization tool and optimality system $\mathtt{(OC)}$ is exact
in this case.

Consider an orthogonal local curvilinear coordinate system coincides
to the principal curvature lines on the $u$-constant surfaces, in
fact the local (intrinsic) tangential coordinate systems. Lets to
denoted by $\textbf{s}(\textbf{x})$ and
$\{\textbf{t}_i(\textbf{x})\}_{i=1}^{d-1}$ respectively the unit
normal and tangent vector which form the basis vectors of this local
curvilinear coordinate system. It is
easy to show that (cf. \cite[][ch. 1]{giga2006see}):%
\begin{equation}%
\label{eq:tniidentitiy}%
   \mathds{1}_d = \textbf{s}(\textbf{x}) + %
   \sum_{i=1}^{d-1} \ \textbf{t}_i(\textbf{x})
\end{equation}
considering \ref{eq:projmat} together with \ref{eq:tniidentitiy}, it
turns out that for every $d$-dimensional vector
$\mathbf{v}(\mathbf{x})$, $P(u(\mathbf{x}))\cdot
\mathbf{v}(\mathbf{x})$ is a $d$-dimensional vector
$\mathbf{w}(\mathbf{x})$ such that $\mathbf{w}(\mathbf{x})$ is
tangential to $u(\textbf{x})$-constant surfaces, i.e.,
$P(u(\mathbf{x}))\cdot \mathbf{v}(\mathbf{x}) \cdot \nabla
u(\mathbf{x}) = 0$. Therefore, at any $\textbf{x}_0$ where $\nabla
\delta u(\textbf{x}_0)$ be orthogonal to the
$u(\textbf{x}_0)$-constant surface, $\delta_u \kappa$ will be equal
to zero. Therefore in a special case when  $\partial \mathcal{D}$ is
defined as  a $u$-constant surface $h_1=j_2=h_3=0$ and so we can
solve the adjoint Poisson equation exactly without any
regularization tool.


\section{Numerical method}%
\label{sec:numerical}%

In this section we briefly mention the numerical method used to
solve the optimality condition $\mathtt{OC}$. The physical domain is
discretized into a uniform Cartesian grid with an $N+1$ grid points
along each spatial dimension. The direct and adjoint systems are
solved by the cell centered finite volume method (FVM). In this way,
we have an $N$ control volumes along each spatial dimension. Each
control volume includes a degree of freedom for direct and adjoint
fields ni addition to a degree of freedom for the control variable.
All degrees of freedom are defined as center of control volumes. The
harmonic averaging is used to compute the diffusion coefficient on
the faces of control volumes. It is worth mentioning that, according
to our experiments and also results reported in
\cite{gersborghansen2006toh}, the cell-centered FVM is free from the
checkerboard-like instability (cf. \cite{bendsoee2004tot}) which is
usually connected to finite element solution of topology
optimization problems.

In the present study, the spectral projected gradient method
developed in \cite{tavakoli2010nonmonotone} is used to solve the
optimality conditions $\mathtt{(OC)}$. To solve the linearized
optimality condition $\mathtt{(OC)}$, the value objective function
together with its gradient is required which needs solution of
direct and adjoint Poisson equations. Due to existence of jump in
diffusion coefficient, the condition number of the coefficient
matrix related to the direct and adjoint problems will be very large
using traditional iterative solvers (in particular for large values
of $k_\beta/k_\alpha$). The multigrid preconditioned conjugate
gradient method (MGCG) \cite{tatebe1993mpc} is used to solve linear
systems of equations in this study. The main benefits of this method
are its excellent performance in addition to its nearly independent
convergence-rate to factor $k_\beta/k_\alpha$.


\section{Results and discussion}%
\label{sec:res}%

The success and performance of the presented approach is studied in
this section through several numerical examples. Some of the input
parameters are as follows. In all of the examples, the design domain
is chosen to be $Q = [-0.5,0.5]^d$. The following strict equality
resource constraint is considered $R_L=R_U=0.5$ and the threshold
for the projection onto the admissible control domain is taken equal
to $1.e-6$. To keep the initial design inside the feasible domain,
we start with $w(\textbf{x})=0.5$. The smoothing width of the
regularization kernel,$\sigma$, is equal to $3\Delta x$, where
$\Delta x$ is the grid spacing. The convergence threshold of MGCG
algorithm is $1.e-20$, except where it is given explicitly. The
jacobi method is used as our multigrid smoother and the number of
smoothing iteration on each grid level (after and before the
prolongation and restriction operations) is one. We consider two
form of definitions for the support of integral in the objective
functional ($\mathcal{D}$). In the first form, $\mathcal{D}$ is
defined by a spatial function and in the second form, $\mathcal{D}$
is taken to be a function of the field solution, for instance, we
first solve the direct problem and then patch a portion of the
spatial domain based on the computed local mean curvature value. In
fact, this type of definition can be considered as an ingredient of
our ultimate goal which is the pointwise control on the second order
derivatives. In this context, the present works plays the role of
the sub-problem solver for an augmented lagrangian method (cf.
\cite{ito2008lma}). The optimization cycle is stopped when the
difference between two consecutive topologies be below $0.1\%$. All
floating point arithmetic is performed to the double precision
accuracy. A personal computer with an AMD 2.41 GHz CPU and 2.5GB RAM
is used as the computing platform.

\subsection{Two dimensional results}%
\label{sec:res:2d}%

The following examples are used to evaluate the presented method in
two spatial dimensions. In these examples, the physical domain is
divided into a $256\times 256$ uniform Cartesian grid. %

%
%
\begin{example}%
\label{ex:2d1}%
$k_\alpha=2$, %
$k_\beta=1$, %
$q=1$, %
$g(\textbf{x})=1$,  %
$u_0(\textbf{x})=0$, %
$a=-1$, %
$b=0$, %
$c=1$, %
$\kappa_0=0$, %
$\mathcal{D} = \{\ \textbf{x} \in [-0.5,0.5]^2 \ \big| \ |x| \le
0.25,\ |y| \le 0.25$\}, where $\textbf{x} = (x, y)^T$.
\end{example}%
%
%
\begin{example}%
\label{ex:2d2}%
this example is like to \ref{ex:2d1}, else $b=1$. %
\end{example}%
%
%
\begin{example}%
\label{ex:2d3}%
this example is like to \ref{ex:2d1}, else the conductivity ratio is
increased to $200$, i.e., $k_\alpha=200$, $k_\beta=1$, the SIMP
power is also increased to $q=5$ avoid intermediate densities.%
\end{example}%
%
%
\begin{example}%
\label{ex:2d4}%
this example is like to \ref{ex:2d1}, else the objective function
domain is defined as a function of the state curvature at the start
of the optimization, i.e., assuming $\kappa_0=-6$ the characteristic
function $\chi(\mathcal{D})$ is defined as:%
\begin{equation}%
\label{eq:2d4}%
    \chi(\mathcal{D}) =
    \left\{%
\begin{array}{rrr}
      1 & \texttt{if} & \kappa(\textbf{x}) \le \kappa_0\\%
      0& \texttt{if} & \kappa(\textbf{x})> \kappa_0%
\end{array}%
\right.%
\nonumber
\end{equation}%
\end{example}%
%
\begin{example}%
\label{ex:2d5}%
this example is like to \ref{ex:2d4}, else $b=1$.%
\end{example}%
\begin{example}%
\label{ex:2d6}%
this example is like to \ref{ex:2d4}, else the direction of the
optimization is reversed, i.e., $a=1$.%
\end{example}%
%
%
\begin{example}%
\label{ex:2d7}%
this example is like to \ref{ex:2d4}, else the conductivity ratio is
increased to $200$, i.e., $k_\alpha=200$, $k_\beta=1$, and the SIMP
power is increased to $5$ too.%
\end{example}%
%
%
\begin{example}%
\label{ex:2d8}%
this example is like to \ref{ex:2d7}, else the direction of the
optimization is reversed, i.e., $a=1$.%
\end{example}%
%
\begin{example}%
\label{ex:2d9}%
this example is like to \ref{ex:2d4}, else the right hand side of
the direct problem is changed to: $g(\textbf{x})=\cos(3\pi x)
\cos(3\pi y)$.
\end{example}%
%
%
\begin{example}%
\label{ex:2d10}%
this example is like to \ref{ex:2d9}, else $b=1$.
\end{example}%
%
%
\begin{example}%
\label{ex:2d11}%
this example is like to \ref{ex:2d9}, else the direction of the
optimization is reversed, i.e., $a=1$; also $q=10$.
\end{example}%
%
\begin{example}%
\label{ex:2d12}%
this example is like to \ref{ex:2d11}, else $b=1$.
\end{example}%
%
\begin{example}%
\label{ex:2d13}%
this example is like to \ref{ex:2d4}, else the asymmetric right hand
side $g(\textbf{x})= \cos(\pi x) \sin(\pi y)$ is applied.
\end{example}%
%
\begin{example}%
\label{ex:2d14}%
this example is like to \ref{ex:2d13}, else $b=1$.
\end{example}%
%

\begin{figure}[ht]%
\begin{center}%
\includegraphics[width=12cm]{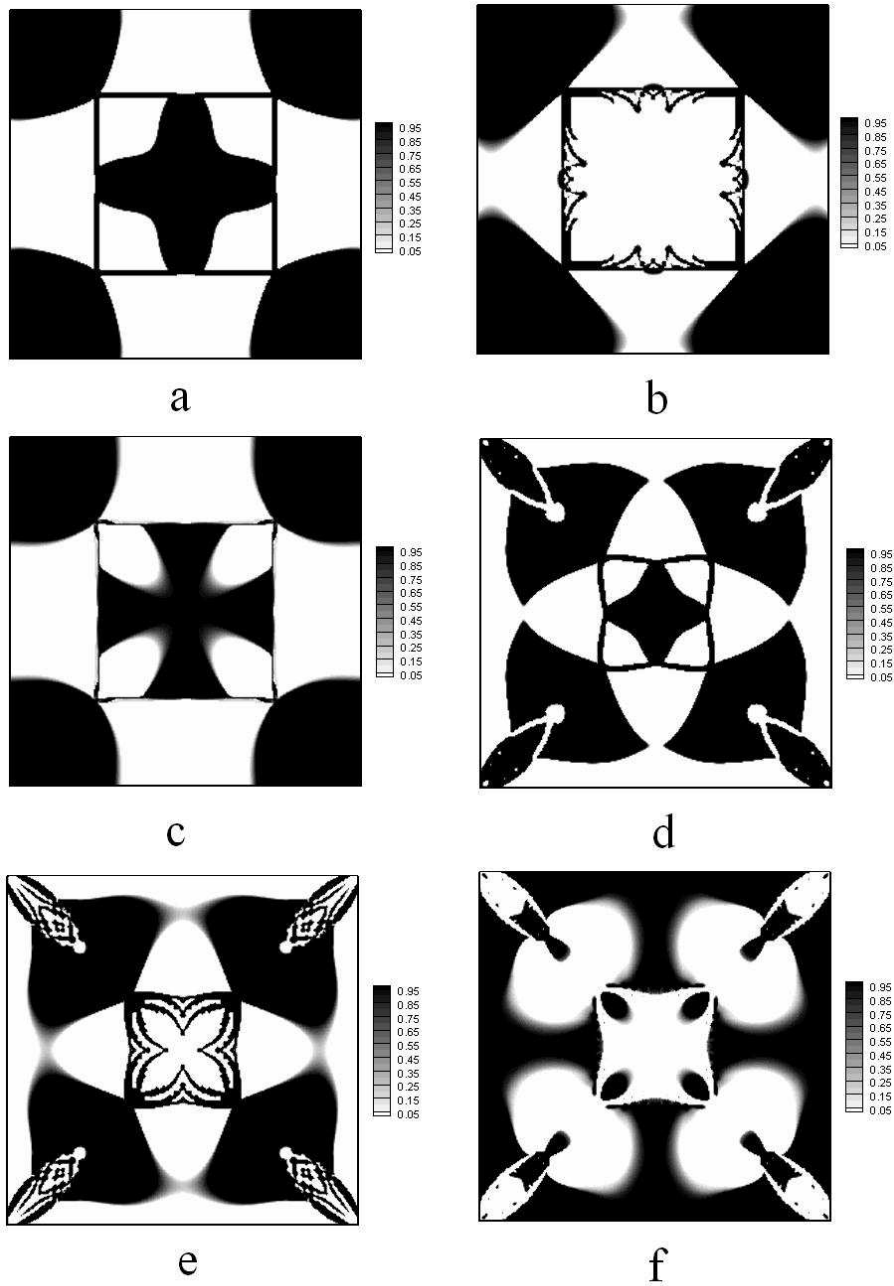}%
\caption{Optimal design in two spatial dimensions:
a-f are related to the final topologies of examples \ref{ex:2d1}-\ref{ex:2d6} respectively.}%
\label{fig:res:2d:ex1_6}%
\end{center}%
\end{figure}

\begin{figure}[ht]%
\begin{center}%
\includegraphics[width=12.cm]{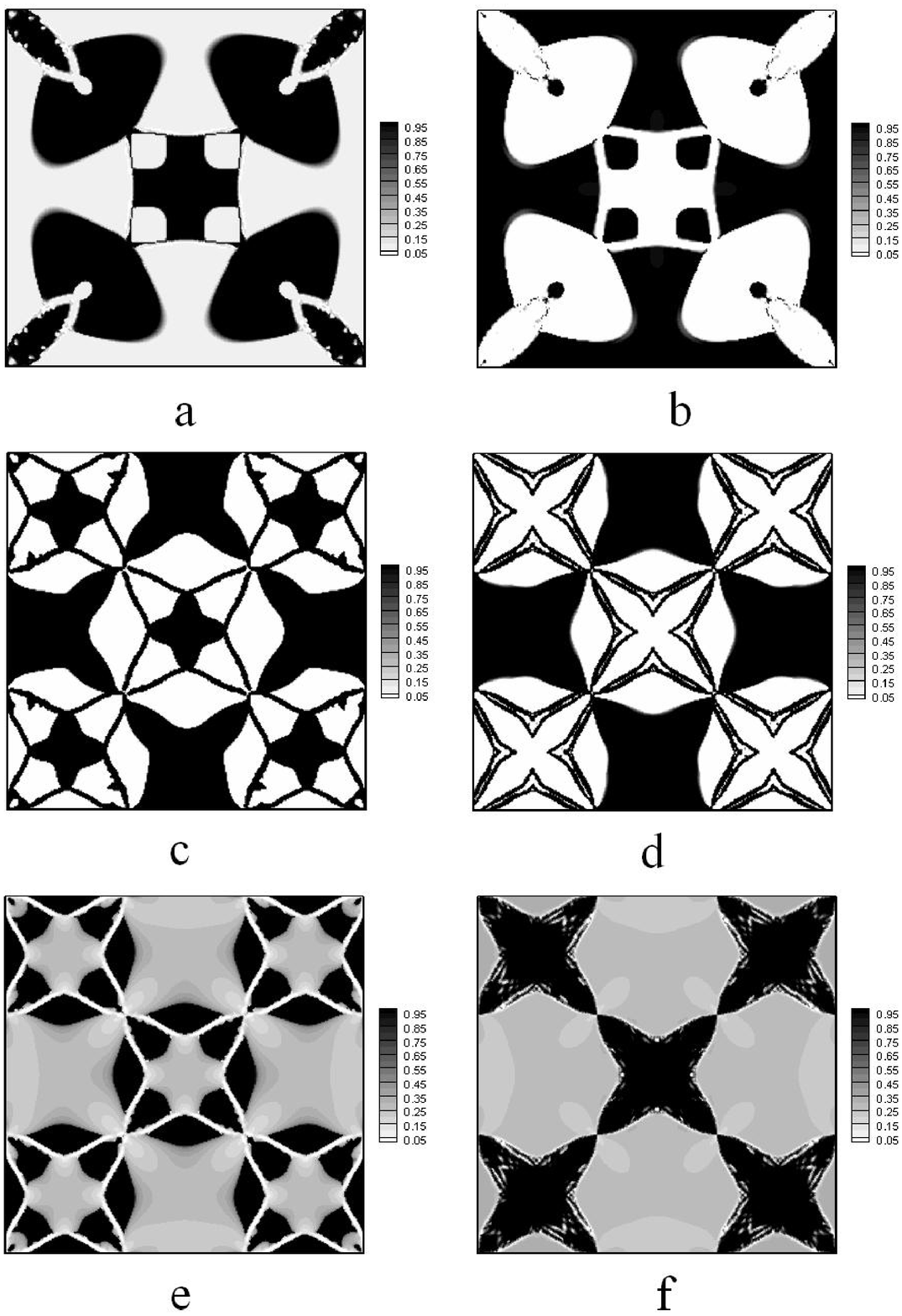}%
\caption{Optimal design in two spatial dimensions:
a-f are related to the final topologies of examples \ref{ex:2d7}-\ref{ex:2d12} respectively.}%
\label{fig:res:2d:ex7_12}%
\end{center}%
\end{figure}

\begin{figure}[ht]%
\begin{center}%
\includegraphics[width=12.cm]{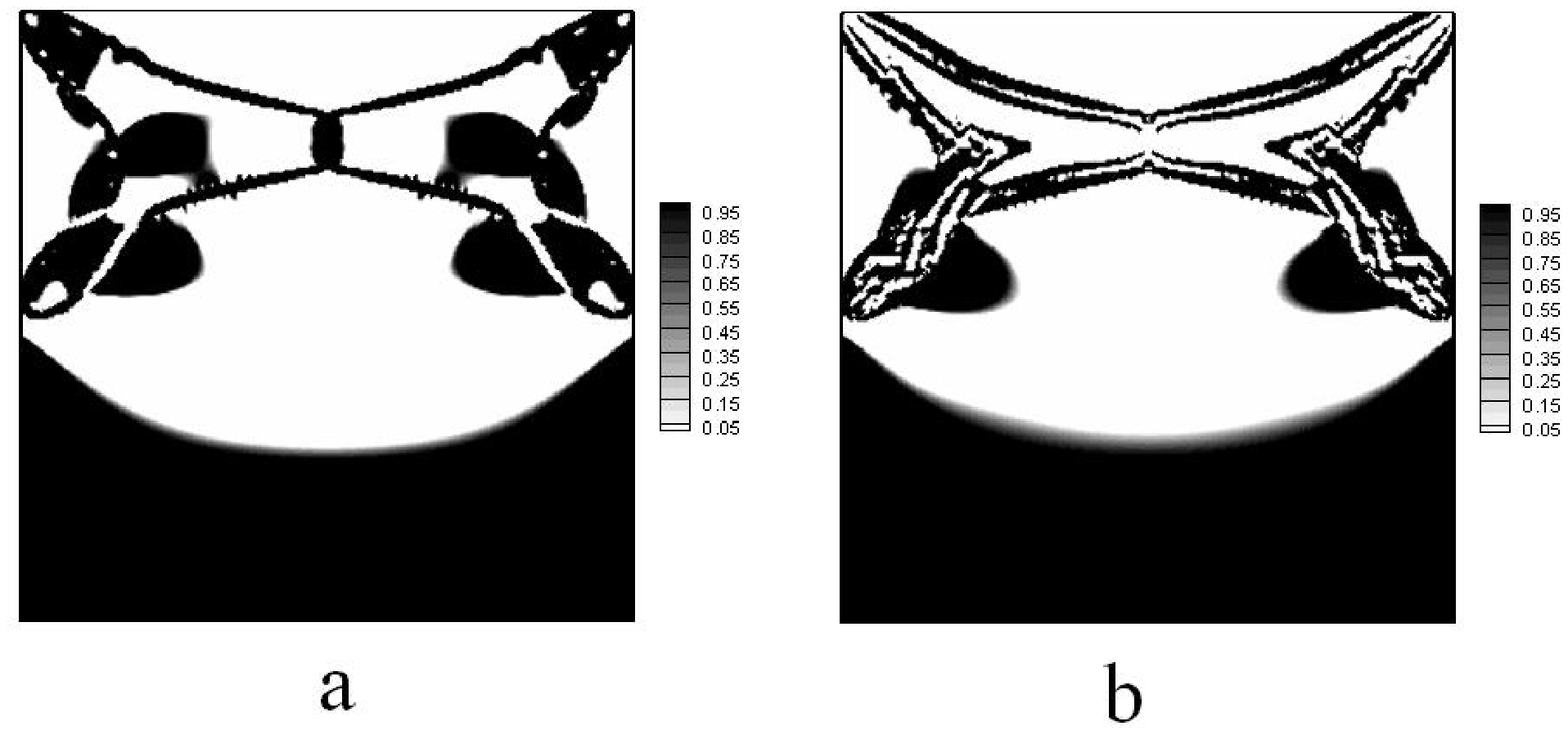}%
\caption{Optimal design in two spatial dimensions:
a and b are related to the final topologies of examples \ref{ex:2d13} and \ref{ex:2d14} respectively.}%
\label{fig:res:2d:ex13_14}%
\end{center}%
\end{figure}

\begin{figure}[ht]%
\begin{center}%
\includegraphics[width=10.5cm]{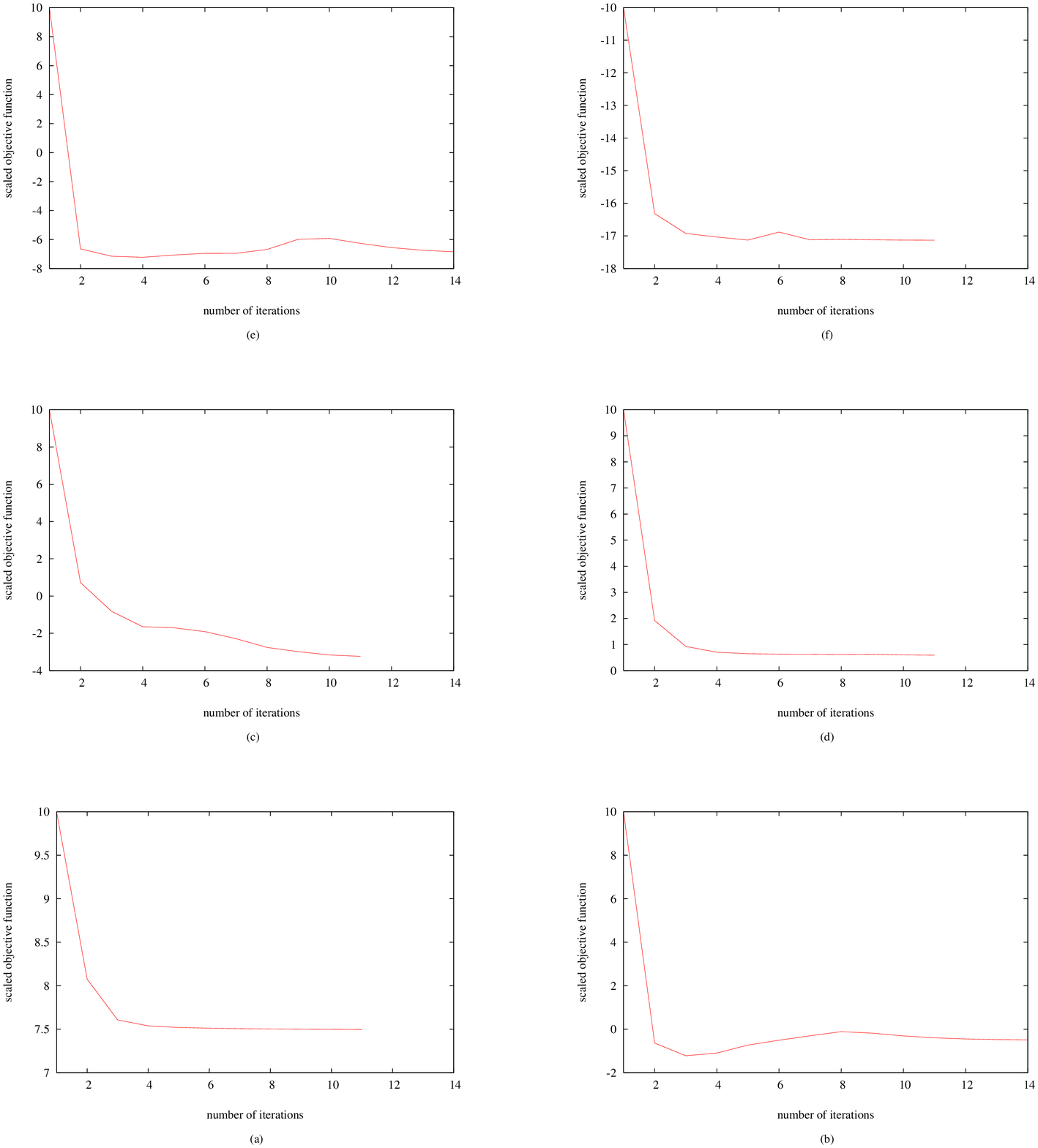}%
\caption{Optimal design in two spatial dimensions:
the objective function history, a-f are related to examples \ref{ex:2d1}-\ref{ex:2d6} respectively.}%
\label{fig:res:2d:o_ex1_6}%
\end{center}%
\end{figure}

\begin{figure}[ht]%
\begin{center}%
\includegraphics[width=10.5cm]{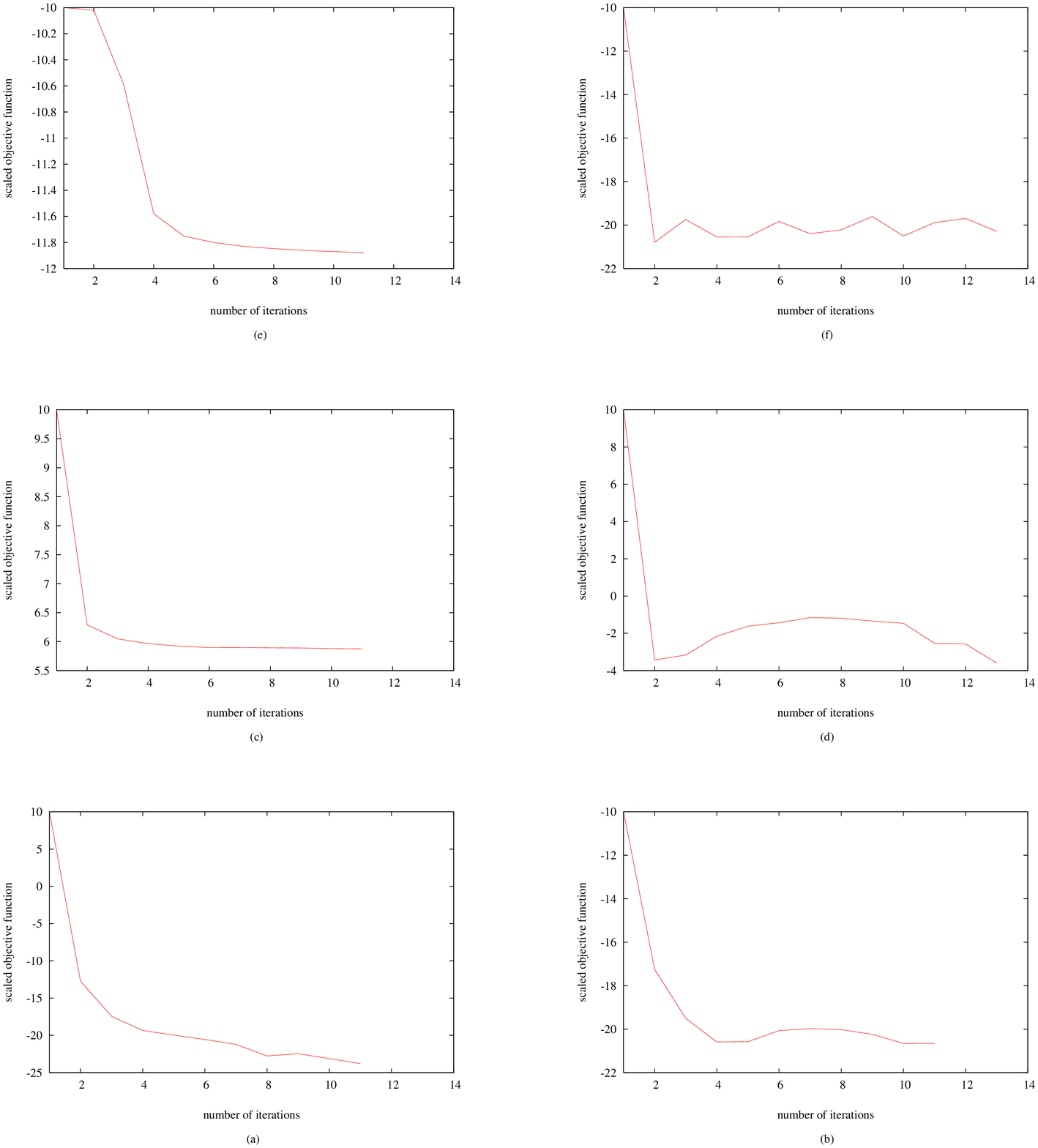}%
\caption{Optimal design in two spatial dimensions: the objective
function history, a-f are related to examples \ref{ex:2d7}-\ref{ex:2d12} respectively.}%
\label{fig:res:2d:o_ex7_12}%
\end{center}%
\end{figure}

\begin{figure}[ht]%
\begin{center}%
\includegraphics[width=10.5cm]{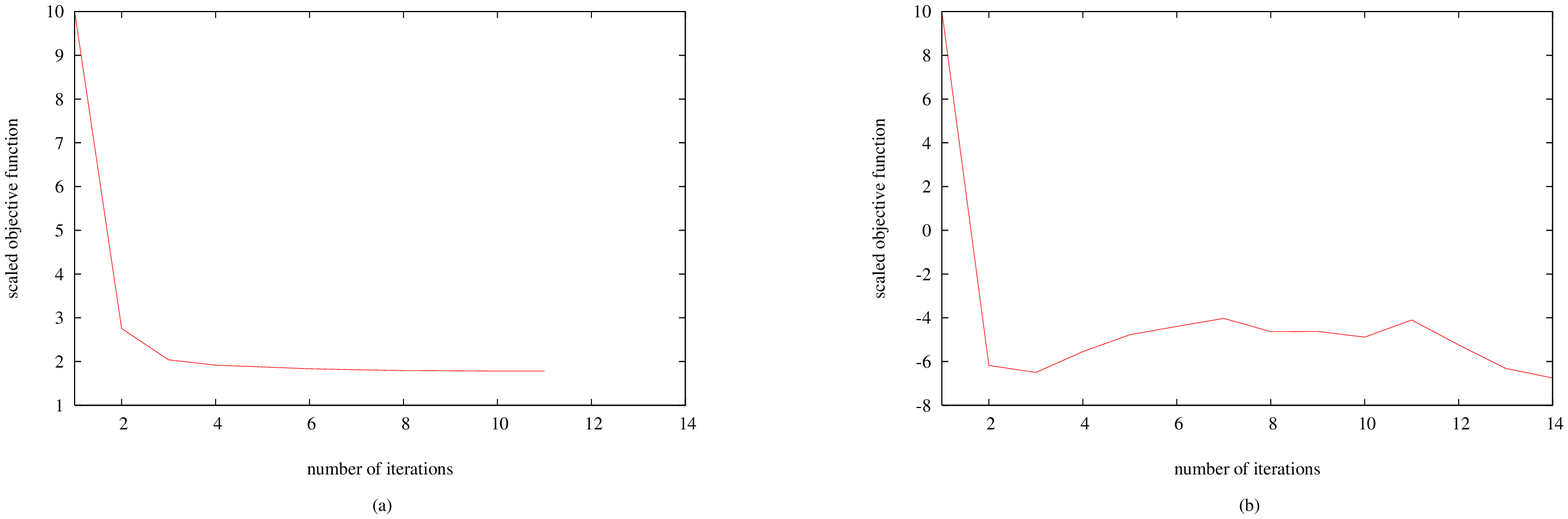}%
\caption{Optimal design in two spatial dimensions: the objective
function history, a and b are related to examples \ref{ex:2d13} and \ref{ex:2d14} respectively.}%
\label{fig:res:2d:o_ex13_14}%
\end{center}%
\end{figure}

The resulted topologies related to examples \ref{ex:2d1}-
\ref{ex:2d14} are shown in Figures
\ref{fig:res:2d:ex1_6}-\ref{fig:res:2d:ex13_14}. Figures
\ref{fig:res:2d:o_ex1_6}-\ref{fig:res:2d:o_ex13_14} show the
variation of the objective functional as the optimization proceeds.
Plots clearly show that, the presented method is enable to
effectively reduce the objective functional in all cases and to move
toward the optimal solution (if there be any). Moreover, in most of
cases a near 0-1 topology is achieved. Therefore, we can conjecture
about the existence of an optimal solution for problem
$\mathtt{(P)}$ defined in this study. As it is clear from the plots,
we do not have essentially a monotonic reduction in the objective
functional value which is an expected result due to the use of a
nonmonotonic line-search in our optimization algorithm (cf.
\cite{tavakoli2010nonmonotone}).

\begin{figure}[ht]%
\begin{center}%
\includegraphics[width=12.37cm]{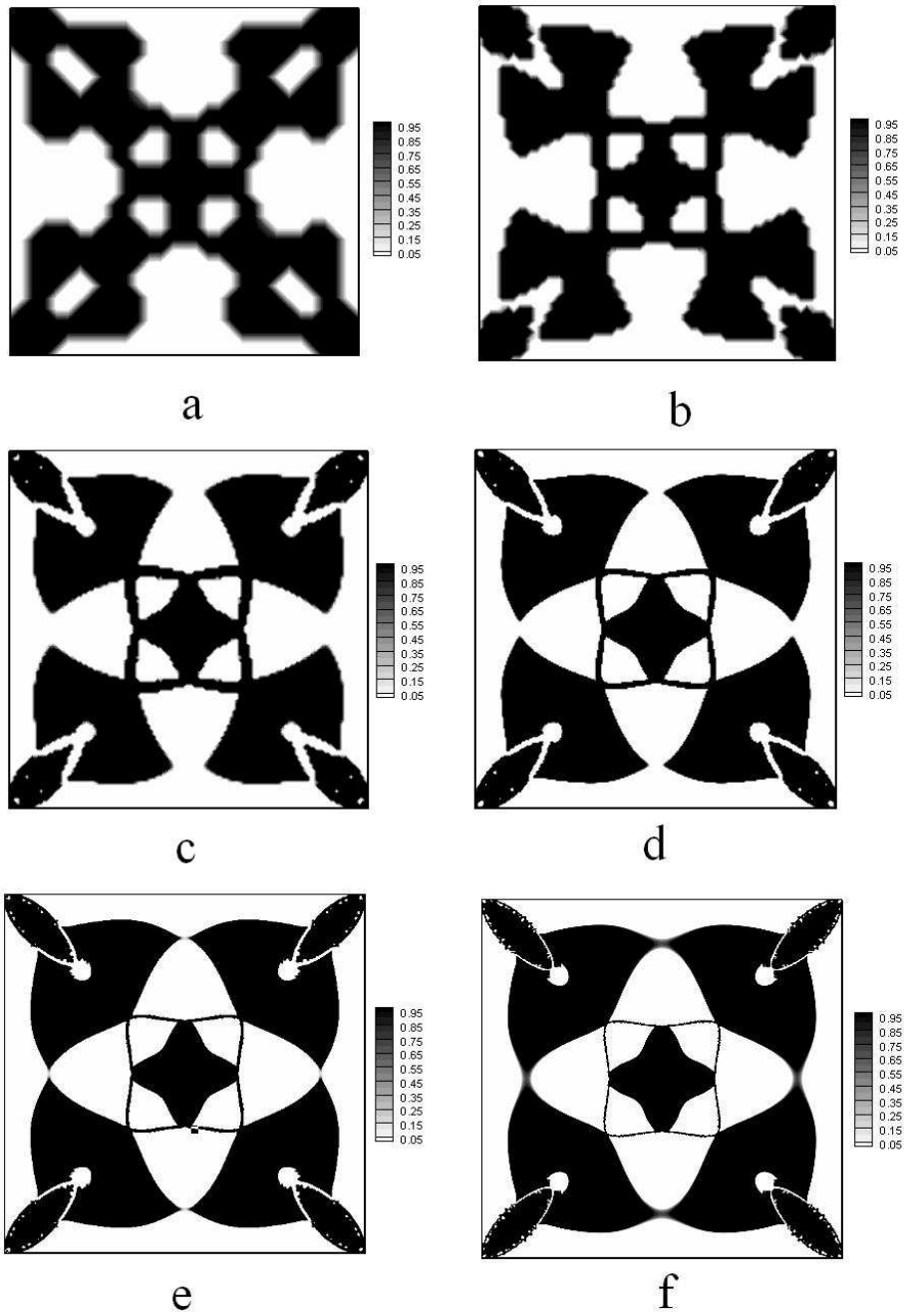}%
\caption{Optimal design in two spatial dimensions: the grid
resolution study, a-f  are related to the final topologies
computed on the grid resolutions $2^n \times 2^n$, $n=5, \cdots, 10$ respectively.}%
\label{fig:res:2d:grd}%
\end{center}%
\end{figure}

To study the stability and convergence of our method with respect to
the grid refinement, example \ref{ex:2d4} is considered with grid
resolutions: $2^n \times 2^n, \ n = 5, \cdots 10$. The resulted
topologies related to this numerical experiment is shown in Figure
\ref{fig:res:2d:grd}. Plot shows that the macroscopic features of
the topology have an excellent convergence. Moreover, the
microscopic features follows the same trend without topological
instability. Therefore, we can conjecture on the well-posedness and
stability of the applied numerical method. Notice that, a
micro-scale topological convergence is generally not expected due to
this well-known fact that by the grid refinement the optimal
solutions tend to form smaller and smaller microstructures (cf.
\cite{allaire2002soh}).

\begin{figure}[ht]%
\begin{center}%
\includegraphics[width=11.cm]{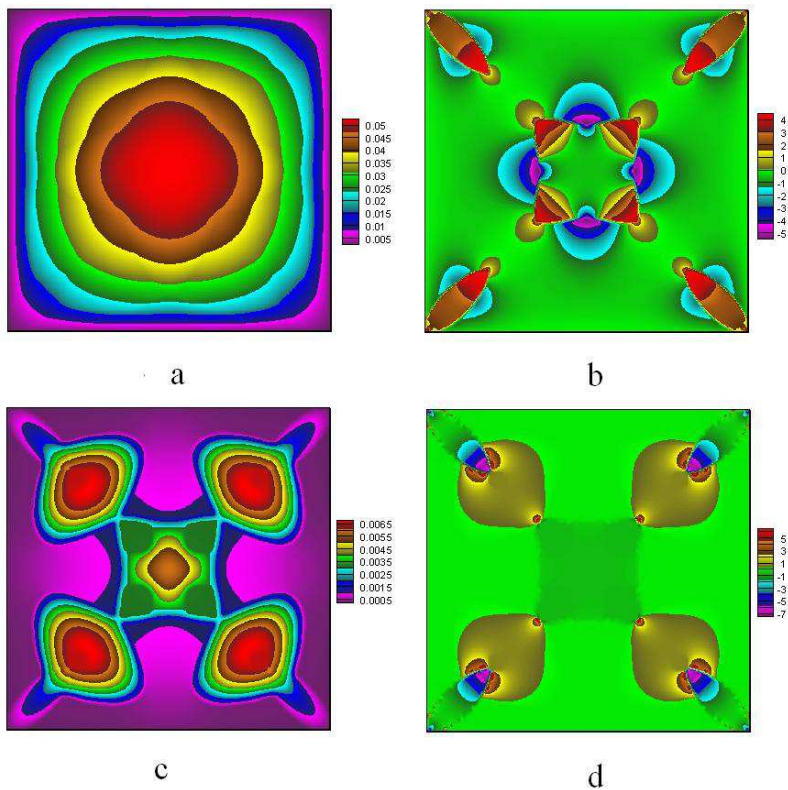}%
\caption{Optimal design in two spatial dimensions: the direct (left)
and adjoint (right) fields at the final design for examples
\ref{ex:2d4} (top) and \ref{ex:2d8} bottom.}%
\label{fig:res:2d:ex48_da}%
\end{center}%
\end{figure}

To make sense about the direct and adjoint fields at optimal
solutions, the direct and adjoint fields at the optimal solution
corresponding to examples \ref{ex:2d4} and
\ref{ex:2d8} are plotted in Figure \ref{fig:res:2d:ex48_da}.%


\clearpage

%


To make sense about the computational cost of the presented
approach,  example \ref{ex:2d4} is considered with different
conductivity ratios $\frac {k_\alpha} {k_\beta} = 2, 10, 100$, and
grid resolutions: $2^n \times 2^n, \ n = 5, \cdots 11$. Notice that
for the conductivity ratio $100$, SIMP power $5$ is considered.

\begin{table}[ht]%
\caption{Performance analysis: results include total optimization
iterations (itero), total functional evaluation (nfunc), total
gradient evaluation (ngard), total MGCG iterations (iterm) and total
computational cost in second (cpu); for different conductivity
ratios ($r_k=k_\alpha/k_\beta$) and different grid resolutions (grid).}%
\label{tab:res:2d:pa}%
\centering %
\begin{tabular}{lllllll}%
\hline\hline%
grid & $r_k$ & itero & nfunc & ngrad & iterm &cpu (sec) \\%
\hline\hline\\%
$(2^6)^2$       & 2&10 & 11 & 11 & 462 & 0.25 \\%
$(2^7)^2$       &2 &10 & 11 & 11 & 484 & 2.15 \\%
$(2^8)^2$       & 2&10 & 11 & 11 & 486 & 10.62 \\%
$(2^9)^2$       & 2&11 & 12 & 12 & 528 & 46.76 \\%
$(2^{10})^2$ & 2 & 11 & 12 & 12 & 542 & 203.1 \\%
$(2^{11})^2$ & 2 &13 & 12 & 12 & 612 & 1180.6 \\\\%
$(2^6)^2$       & 10&10 & 11 & 11 & 792 & 0.36 \\%
$(2^7)^2$       & 10&10 & 11 & 11 & 748 & 3.16 \\%
$(2^8)^2$       & 10&10 & 11 & 11 & 836 & 16.25 \\%
$(2^9)^2$       & 10&11 & 12 & 12 & 990 & 78.30 \\%
$(2^{10})^2$ & 10 &11& 12 & 12 & 1020 & 382.4 \\%
$(2^{11})^2$ & 10 & 14&15 & 15 & 1172 & 1905.1 \\\\%
$(2^6)^2$       & 100 &10 & 11 & 11 & 1430 & 0.56 \\%
$(2^7)^2$       & 100 &10 & 11 & 11 & 1716 & 6.49 \\%
$(2^8)^2$       & 100 &10 & 11 & 11 & 1958 & 35.05 \\%
$(2^9)^2$       & 100 &12 & 13 & 13 & 3014 & 220.7 \\%
$(2^{10})^2$ & 100 &12 & 13 & 13 & 4940 & 1573.5 \\%
$(2^{11})^2$ & 100 &14 & 15 & 15 & 6582 & 9941.1\\%
\hline\hline%
\end{tabular}%
\end{table}%

Results of this numerical experiment are shown in Table
\ref{tab:res:2d:pa}. Table illustrates that the number of
optimization cycles has a little dependency on the grid resolution
and the conductivity ratio. The number of MGCG iterations is
increased by increasing the conductivity ratio, $r_k$. But for a
fixed $r_k$, the number of MGCG iterations is almost independent
from the grid resolution, which show the reliability of this solver
for large-scale problems. The reported computational cost
illustrates the excellent performance of the applied numerical
strategy in this study. To the best of our knowledge such a large
number of design parameters ($2^{22}$) is not reported in literature
of topology optimization so far.


\subsection{Three dimensional results}%
\label{sec:res:3d}%

In this section we show the success of the presented method in three
spatial dimensions. For this purpose the following example is
considered on grid resolutions: $2^n \times 2^n \times 2^n$, $n=5,
6, 7$.

%
%
\begin{example}%
\label{ex:3d1}%
$k_\alpha=2$, %
 $k_\beta=1$,  %
$q=1$, %
$g(\textbf{x})=1$, %
$u_0(\textbf{x})=0$, %
$a=-1$, %
$b=0$, %
$c=1$, %
and assuming $\kappa_0=-3$, the characteristic function
$\chi(\mathcal{D})$ is defined as%
\begin{equation}%
\label{eq:2d4}%
    \chi(\mathcal{D}) =
    \left\{%
\begin{array}{rrr}
      1 & \texttt{if} & \kappa_\varepsilon(\textbf{x}) \le \kappa_0\\%
      0& \texttt{if} & \kappa_\varepsilon(\textbf{x})> \kappa_0%
\end{array}%
\right.%
\nonumber
\end{equation}%
\end{example}%

\begin{figure}[ht]%
\begin{center}%
\includegraphics[width=10cm]{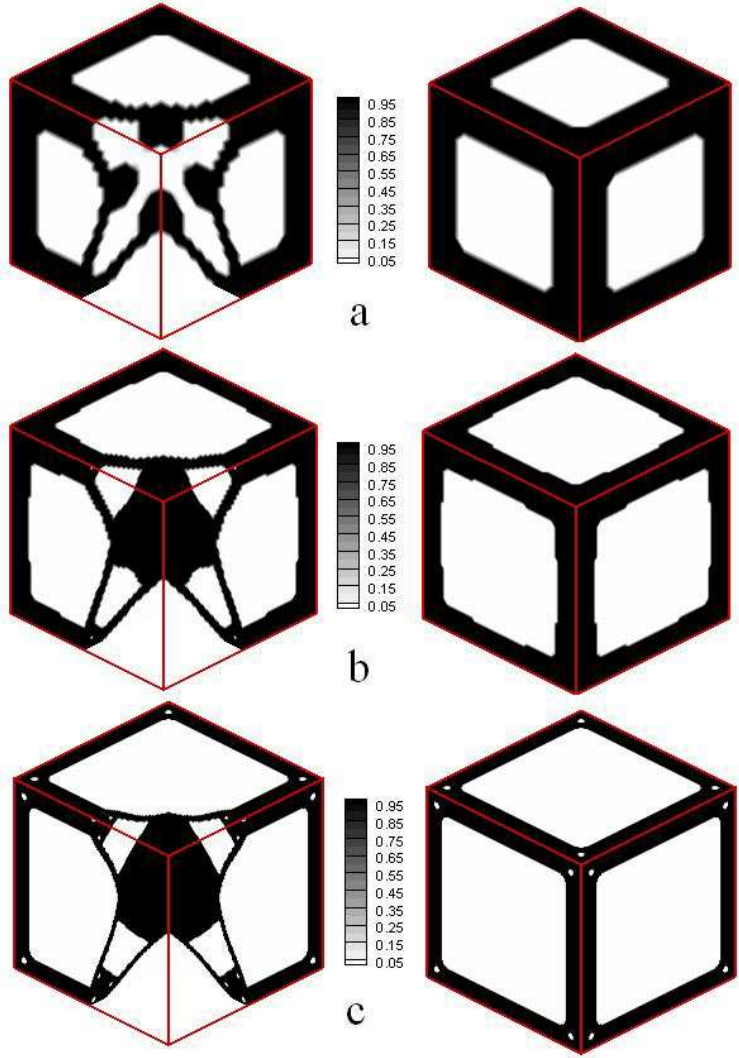}%
\caption{Optimal design in three spatial dimensions: a-c are related
to the final topologies computed on the grid resolutions
$2^n \times 2^n \times 2^n$, $n=5, 6, 7$ respectively.}%
\label{fig:res:3d:3d}%
\end{center}%
\end{figure}

\begin{figure}[ht]%
\begin{center}%
\includegraphics[width=13cm]{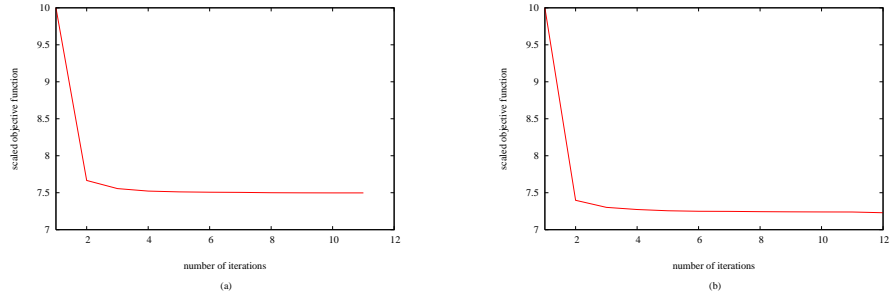}%
\caption{Optimal design in three spatial dimensions: the objective
function history, a and b are related to grid resolutions
$2^6 \times 2^6 \times 2^6$ and $2^7 \times 2^7 \times 2^7$ respectively.}%
\label{fig:res:3d:o_3d}%
\end{center}%
\end{figure}

Figures \ref{fig:res:3d:3d} and \ref{fig:res:3d:o_3d}   show the
final topologies and objective function history corresponding to
three dimensional numerical examples. In all cases the number of
optimization iterations was below 13 and the number of objective
functional and its gradient evaluation was equal. Plots illustrate
the convergence and stability of the presented method in three
spatial dimensions too.

\clearpage

%


\section{Summary}%
\label{sec:remarks}%

The idea of posing control on the second order derivatives of the
PDEs solutions is introduced in this study. To be specific, a
topology optimization problem corresponding to the Poisson equation
in which the objective functional is  a nonlinear function of first
and second order field derivatives is considered. Introducing some
functional regularization tools, the (approximate) first order
necessary optimality conditions is derived and is numerical solved
using an appropriate gradient descent method. Numerical results are
provided for a variety of test cases in two and three spatial
dimensions. The Numerical results in this study, confirms the
stability, convergence and efficiency of the presented approach. In
fact, the numerical results answer to our quest in this study; the
possibility of posing some degrees of control on the second order
derivatives.

\section*{Acknowledgment}

The author would like to thanks Wolfgang Bangerth for his
constructive comments in the course of this research, also to
Hongchao Zhang for his comments on the implementation of projected
gradient method. The implementation of MGCG method in this study is
inspired from the implementation by Osamu Tatebe which is available
from his web. The lsmlib library of Kevin T. Chu is applied to
compute the distance filed and field extension in this study. The
assistance of Kevin T. Chu in this regard is also appreciated.


\clearpage

\bibliographystyle{plain} 
\bibliography{biblio}%

\end{document}